\documentclass[12pt]{amsart}

\input{epsf}
\usepackage{amssymb,latexsym}

\font\smallit=cmti10
\font\smalltt=cmtt10
\font\smallrm=cmr9

\newtheorem{theorem}{Theorem}[section]

\newtheorem{lemma}[theorem]{Lemma}

\begin{document}
\pagenumbering{arabic}
\pagestyle{headings}
\def\sof{\hfill\rule{2mm}{2mm}}
\def\NN{\mathbb{N}}
\def\mn{\mbox{-}}

\begin{center}
{\bf ENUMERATION OF 3-LETTER PATTERNS IN COMPOSITIONS}
\vskip 20pt
{\bf Silvia Heubach}\\
{\smallit Department of Mathematics, California State University Los Angeles,
Los Angeles, CA 90032, USA}\\
{\tt sheubac@calstatela.edu}\\ 
\vskip 10pt
{\bf Toufik Mansour}\\
{\smallit Department of Mathematics, Haifa University, 31905 Haifa, Israel}\\
{\tt toufik@math.haifa.ac.il}\\ 
\end{center}
\vskip 30pt
\centerline{\smallit Received: , Accepted: , Published: }
\vskip 30pt

\centerline{\bf Abstract}

\noindent
Let $A$ be any set of positive integers and $n\in\NN$. A {\em composition of
$n$ with parts in $A$} is an ordered collection of one
or more elements in $A$ whose sum is $n$. We derive generating functions for the number of compositions of $n$ with $m$ parts in $A$ that have $r$ occurrences of 3-letter patterns formed by two (adjacent) instances of {\it levels}, {\it rises} and {\it drops}. We also derive asymptotics for the number of compositions of $n$ that avoid a given pattern. Finally, we obtain the generating function for the number of $k$-ary words of length $m$ which contain a prescribed number of occurrences of a given pattern as a special case of our results.



\pagestyle{myheadings}
\markright{\smalltt INTEGERS: \smallrm ELECTRONIC JOURNAL OF COMBINATORIAL NUMBER THEORY \smalltt x (200x), \#Axx\hfill}

\thispagestyle{empty} 
\baselineskip=15pt 
\vskip 30pt

\vskip 30pt
\section{Introduction}
A {\em composition} $\sigma=\sigma_1\sigma_2\ldots\sigma_m$ of $n\in\NN$
is an ordered collection of one or more positive
integers whose sum is $n$. The number of summands, namely $m$, is called
the number of {\em parts} of the composition.
We will look at {\em compositions of $n$ with parts in $A$}, i.e.,
compositions whose parts are restricted to be from a set $A \subseteq
\NN$. Our aim is to count the number of compositions of $n$ with parts in
$A$ which contain a $3$-letter pattern $\tau$ exactly $r$ times. This extends work on the statistics {\em rises} (a summand followed by a larger summand), {\em levels} (a summand
followed by itself), and {\em drops} or {\em falls} (a summand followed by a smaller
summand) in all compositions of $n$ whose parts are in a given set $A$.

Several authors have studied compositions and the statistics rises, levels and falls from different viewpoints.
Alladi and Hoggatt ~\cite{AllHog1975} studied compositions with parts in the set $A=\{1,2\}$ in conjunction with the
Fibonacci sequence. Chinn, Grimaldi and Heubach ~\cite{ChiGriHeu2003, ChiHeu2003a, ChiHeu2003, Gri2000, Gri2001} have generalized to different sets $A$ and have counted the number
of compositions, rises, levels and drops, looking for connections to known sequences (which is only possible when considering a specific set $A$).

Carlitz and several co-authors
(\cite{Car1973},\cite{Car1977},\cite{CarSco1972},\cite{CarScoVau1976},\cite{CarVau1974}) studied rises,
levels and falls in compositions on the set $[n]=\{1,2,\ldots,n\}$ as an
extension of the study of these statistics or patterns in permutations. (Unlike Alladi and Hoggatt, Carlitz et al.
included an additional rise at the beginning and an additional fall at the end of each composition, except in~\cite{CarScoVau1976}.) These authors extended enumeration questions for permutations to compositions by also considering the {\em specification} of a composition (a list of counts for the occurrences of each integer) as well as the statistic levels. More recently, Rawlings~\cite{Raw2000} enumerated compositions according to weak rises and falls
(equality allowed) in connection with restricted words by adjacencies. He also introduced the notion of {\em ascent
variation}, (the sum of the increases of the rises within a composition), which is motivated by a connection to the
perimeter of directed vertically convex polyominoes. Furthermore, Heubach and Mansour~\cite{HeuMan2004} developed a
general framework which gives the generating function for the number of rises, levels and falls  for any ordered subset $A$ of $\NN$.

``Closest'' to permutations are those compositions that do not have levels, which were called  {\em waves} by Carlitz
and his co-authors.  These compositions have also been referred to as {\em Smirnov sequences}~(see for
example~\cite{GouJac1983}) and  {\em Carlitz
compositions}~\cite{KnoPro1998}. One special case of Carlitz  compositions are up-down sequences, in which rises and falls alternate, which were studied by Carlitz and
Scoville~\cite{CarSco1972} and Carlitz~\cite{Car1973}. In addition, the problem of enumerating Carlitz compositions
according to rises and falls reduces to the Simon Newcomb problem (\cite{Car1972}, \cite{DilRos1969}) when the  number of falls is disregarded.    Another extension of questions first studied for permutations are the  generating
functions for
the number of compositions according to specification, rises,
falls and maxima by Carlitz and Vaughan~\cite{CarVau1974}, and enumeration of pairs of sequences according to rises,
levels and falls \cite{CarScoVau1976}. Recently, Gho, Hitczenko, Louchard and Prodinger have studied distinctness and other characteristics of compositions and Carlitz compositions using a probabilistic approach~\cite{GohHit2002, HitLou2001, LouPro2002}.

A widely studied topic for permutations is pattern avoidance (see~\cite{Bon2004} and references therein for an overview). We follow the motivation of extending enumeration questions from permutations to compositions by expressing rises, levels and drops as 2-letter patterns. For example, the pattern 11 corresponds to any
occurrence of $a_ia_i$ with $a_i \in A$, and thus any occurrence of the pattern 11 corresponds to the
occurrence of a level. Likewise, the pattern 12 corresponds to $a_ia_j$ with $a_i < a_j$, i.e., a rise, and the pattern 21 corresponds to a
fall or drop.

In this paper, we generalize to 3-letter patterns. For example, the pattern
111 corresponds to any occurrence of $a_ia_ia_i$ with $a_i \in A$,
and therefore corresponds to a level followed by a level; we will
refer to this statistic by the shorthand level+level. Likewise, we
can define patterns for all the combinations of rises, levels and
drops. To illustrate this idea, we look at the composition
1413364. In terms of rises, levels and drops, this composition is
represented by rise+drop+rise+level+rise+drop. It contains the
patterns (from left to right) 121 (rise+drop), 312 (drop+rise),
122 (rise+level), 112 (level+rise) and 132 (rise+drop).
\newline
\newline Due to symmetry, each rise is matched by a drop (for each
composition that is not symmetric, there is a composition whose
parts are in reverse order, and for symmetric compositions, the
rises and drops are matched within the composition). Thus, the
statistic rise+level, 122, occurs as often as the statistic
level+drop, 221. Table~\ref{patt} lists the statistics to be considered,
and their corresponding patterns. Note that the statistic rise+drop is represented by
three different patterns, which take into account the relative
size of the actual summands. For example, the pattern 121 indicates the
occurrence of $a_1a_2a_1$ where $a_1<a_2$, whereas 132 indicates
the occurrence of $a_1a_2a_3$, where $a_1<a_3<a_2$.

\begin{table}
  \begin{center}
\begin{tabular}{|c|c||c|c|} \hline
Statistic& Pattern&Statistic& Pattern\\ \hline
level+level&111&rise+rise&123\\
level+rise&112 &rise+drop=peak&121+132+231\\
level+drop&221& drop+rise=valley&212+213+312\\ \hline
\end{tabular}
\end{center}
  \caption{Statistics and their associated patterns}\label{patt}
\end{table}

In Section~\ref{gfs} we derive the generating functions for the number of compositions of $n$ with $m$ parts in $A$ which contain a given pattern $\tau$ exactly $r$ times for each of the patterns listed above. In Section~\ref{asymps} we use tools from complex analysis to derive the exact asymptotics for the number of compositions of $n$ that avoid a given pattern $\tau$. Finally, in  Section~\ref{words}, we apply our results to words, i.e., elements of $[k]^n$, where $[k]$ is a (totally ordered) alphabet on $k$ letters.  We obtain previous results in~\cite{BurMan2003, BurMan2003a}, and obtain new results for the patterns peak and valley. This application shows that compositions form a larger class of combinatorial objects, containing words as a subclass.

\vskip 30pt
\section{Compositions of $n$ with parts in $A$} \label{gfs}
Let $A=\{a_1,a_2,a_3,\ldots,a_d\}$ or $A=\{a_1,a_2,a_3,\ldots\}$,
where $a_1<a_2<\ldots.$ are positive integers. We will refer to such
a set as an {\em ordered subset} of $\NN$. In the theorems and
proofs, we will treat the two cases together if possible, and will
note if the case $|A|=\infty$ requires additional steps.

Let $C_{\tau}(n,r)$ (respectively $C_{\tau}(j;n,r)$) denote the number
of compositions of $n$ with parts in $A$ (respectively with $j$
parts in $A$) which contain the pattern $\tau$ exactly $r$ times.
The corresponding generating functions are given by
$$C_\tau(x,y)=\sum_{n,r\geq0}C_\tau(n,r)x^ny^r$$
and
$$
C_\tau(x,y,z)=\sum_{n,r,j\geq0}C_\tau(j;n,r)x^ny^rz^j=\sum_{j\geq0}C_\tau(j;x,y)z^j.$$

More generally, let $C_{\tau}(\sigma_1\ldots\sigma_\ell|n,r)$
(respectively $C_{\tau}(\sigma_1\ldots\sigma_\ell|j;n,r)$) be the
number of compositions of $n$ with parts in $A$ (respectively with
$j$ parts in $A$) which contain $\tau$ exactly $r$ times and whose first $\ell$ parts are
$\sigma_1,\ldots,\sigma_\ell$. The corresponding generating
functions are given by
$$C_\tau(\sigma_1\ldots\sigma_\ell|x,y)=
\sum_{n,r\geq0}C_\tau(\sigma_1\ldots\sigma_\ell|n,r)x^ny^r$$
and
$$C_\tau(\sigma_1\ldots\sigma_\ell|x,y,z)=
\sum_{n,r,j\geq0}C_\tau(\sigma_1\ldots\sigma_\ell|j;n,r)x^ny^rz^j=
\sum_{j\geq0}C_\tau(\sigma_1\ldots\sigma_\ell|j;x,y)z^j.$$ 
The
initial conditions are $C_{\tau}(j;x,y)=0$ for $j<0$,
$C_{\tau}(0;x,y)=1$, and 
$C_{\tau}(\sigma_1, \dots \sigma_{l}|j;x,y)=0$ for $j \leq l-1$. In addition,
\begin{equation}\label{eqfact1}
C_{\tau}(x,y,z)=1+\sum_{a\in A}C_{\tau}(a|x,y,z).
\end{equation}

In this section we study the generating functions
$C_{\tau}(x,y,z)$ for different patterns
$\tau=\tau_1\tau_2\tau_3$. To find an explicit expression for
$C_{\tau}(x,y,z)$, we derive recursive equations using a variety of strategies for the different patterns. For patterns that contain levels, namely $111$, $112$, and $221$, the recursion is in terms of
$C_{\tau}(\sigma_1|j;x,y)$, $C_{\tau}(\sigma_1\sigma_2|j;x,y)$,
$C_{\tau}(\sigma_1\sigma_2\sigma_3|j;x,y)$ and $C_{\tau}(x,y,z)$,
which is usually pretty straightforward. However, solving the
resulting system of equations can be difficult, as for
example in the case of the pattern 112. We will describe the
derivation of the recursive equations in detail for the pattern
111.

For the pattern $123$, we break the composition of $n$ into pieces, some of which have parts in $A$, and others that contain only parts larger than the part under consideration. We then define a second generating function which will play a major role in the recursive equation for $C_{\tau}(x,y,z)$. In the case of peaks (valleys), we split the composition into parts according to where the largest (smallest) part occurs, and derive a recursion that will lead to a continued fraction expansion for the generating function $C_{\tau}(x,y,z)$.

\subsection{The pattern $111$ (the statistic level+level)}
In the following theorem we present the generating function for
the number of compositions of $n$ with $j$ parts in $A$ that
contain the pattern $111$ exactly $r$ times.

\begin{theorem}\label{th111} Let $A$ be any ordered subset of $\NN$. Then
$$C_{111}(x,y,z)=\frac{1}{1-\sum_{a\in
A}\frac{x^az(1+(1-y)x^az)}{1+x^az(1+x^az)(1-y)}}.$$
\end{theorem}
{\it Proof.}
The pattern 111 occurs when $a \in A$ occurs three times in a row. Thus, for fixed $a \in A$ and $j \ge 2$
\begin{equation}\label{eq111a}
\begin{array}{l}
C_{111}(a|j;x,y)=C_{111}(aa|j;x,y)+\sum_{b\in A,b\neq a}C_{111}(ab|j;x,y)\\
\qquad=C_{111}(aa|j;x,y)+x^aC_{111}(j-1;x,y)-x^aC_{111}(a|j-1;x,y).
\end{array}\end{equation}
Note that the factor of $x^a$ reflects the fact that we are
looking at compositions of $n-a$. We now apply a similar argument to $C_{111}(aa|j;x,y)$ to obtain  for $j \ge 3$
\begin{equation}\label{eq111b}
\begin{array}{ll}
C_{111}(aa|j;x,y)=C_{111}(aaa|j;x,y)+\sum_{b\in A,b\neq a}C_{111}(aab|j;x,y)\\
\qquad=x^a\,yC_{111}(aa|j-1;x,y)+x^{2a}\sum_{b\in A,b\neq a}C_{111}(b|j-2;x,y)\\
\qquad=x^a\,yC_{111}(aa|j-1;x,y)+x^{2a}C_{111}(j-2;x,y)-x^{2a}C_{111}(a|j-2;x,y).
\end{array}\end{equation}
Multiplying (\ref{eq111a}) and (\ref{eq111b}) by $z^j$, summing over
all $j\geq1$, taking into account that the recurrences hold for $j \ge 2$ and $j \ge 3$, 
and solving the resulting system of two equations
for $C_{111}(a|x,y,z)$, we get that
$$C_{111}(a|x,y,z)=\frac{x^az(1+(1-y)x^az)}{1+x^az(1+x^az)(1-y)}
C_{111}(x,y,z).$$ Summing over all $a\in A$ and using
Eq.~(\ref{eqfact1}) for $\tau=111$, we get the desired result.
\sof

Applying Theorem~\ref{th111} to $A=\NN$ with $a_i=i$ for $i\ge1$, we get that the
generating function for the number of compositions of $n$ with
parts in $\NN$ which avoid the pattern 111 is given by
$$C_{111}(x,0,1)=\frac{1}{1-\sum_{i\geq1}\frac{x^i(1+x^i)}{1+x^i(1+x^i)}},$$
and the values of the corresponding sequence are 1, 1, 2, 3, 7,
13, 24, 46, 89, 170, 324, 618, 1183, 2260, 4318, 8249, 15765,
30123, 57556, 109973, 210137, 401525, 767216, 1465963, 2801115,
5352275 for $n=0,1,\ldots,25$. Note that
compositions that avoid the pattern 111 have only isolated levels.

{\bf Remark: } We note that {\em Carlitz compositions} of $n$, introduced in~\cite{Car1976}, are those compositions of $n$ in which no adjacent parts are the same. Thus, Carlitz compositions are precisely those compositions that avoid levels, or equivalently, avoid the (2-letter) pattern 11. One possible generalization of Carlitz compositions is to define {\it $\ell$-Carlitz compositions of $n$} to be those compositions of $n$ that avoid  $\ell$ consecutive levels, or in terms of pattern avoidance, avoid the pattern $11\ldots11$ consisting of $\ell+1$ 1's.

\subsection{The patterns $112$  and $221$ (the statistics level+rise and level+drop)}
In the following theorem we present the generating functions for
the number of compositions of $n$ with $j$ parts in $A$ which
contain the patterns $112$ and $221$, respectively, exactly $r$ times.

\begin{theorem}\label{th112}
Let $A$ be any ordered subset of $\NN$. Then 
 $$C_{112}(x,y,z)= \frac{1}{1-\sum_{j=1}^d
\left( x^{a_j}z\prod_{i=1}^{j-1}(1-(1-y)x^{2a_i}z^2)\right)}$$
and
$$C_{221}(x,y,z)=\frac{1}{1-\sum_{j=1}^d\left(x^{a_j}z\prod_{i=j+1}^d(1-(1-y)x^{2a_i}z^2)\right)}.$$
\end{theorem}

{\it Proof.}  To derive the generating function $C_{112}(x,y,z)$ we use arguments similar to those in the proof of Theorem~\ref{th111} to obtain for every $a \in A$
\begin{align} \label{c112eq}
C_{112}(a|x,y,z)&=\frac{x^{2a}z^2}{1-x^{2a}z^2}+\frac{x^{2a}z^2}{1-x^{2a}z^2}\sum_{b\in A, b <
a}C_{112}(b|x,y,x) \nonumber \\
&\phantom{=}+\frac{x^{2a}z^2y}{1-x^{2a}z^2}\sum_{b\in A,b > a}C_{112}(b|x,y,z)+\frac{x^az}{1+x^az}C_{112}(x,y,z).
\end{align}
Let's now assume that $A$ is finite, i.e. $A=\{a_1,\ldots,a_d\}$. Setting $x_{0}=C_{112}(x,y,z)$, $x_i=C_{112}(a_i|x,y,z)$, $\alpha_i=\frac{x^{2a_i}z^2}{1-x^{2a_i}z^2}$, and
$\beta_i=\frac{x^{a_i}z}{1+x^{a_i}z}$, the above equation is of the form
$$x_i-\alpha_i\sum_{j<i}x_j-\alpha_iy\sum_{j>i}x_j-\beta_ix_0=\alpha_i \quad \mbox{for} \quad i=1,\ldots,d.$$
Together with Eq.~(\ref{eqfact1}) for $\tau=112$ this results in the following system of equations:
{\small \begin{equation*} 
\left(\begin{array}{ccccccc} \label{sy112a}
-\beta_1&1&-\alpha_1y&-\alpha_1y&\cdots&-\alpha_1y&-\alpha_1y\\
-\beta_2&-\alpha_2&1&-\alpha_2y&\cdots&-\alpha_2y&-\alpha_2y\\
-\beta_3&-\alpha_3&-\alpha_3&1&\cdots&-\alpha_3y&-\alpha_3y\\
-\beta_4&-\alpha_4&-\alpha_4&-\alpha_4&\cdots&-\alpha_4y&-\alpha_4y\\
\vdots&&&\vdots&&&\vdots\\
-\beta_{d-1}&-\alpha_{d-1}&-\alpha_{d-1}&-\alpha_{d-1}&\cdots&1&-\alpha_{d-1}y\\
-\beta_d&-\alpha_d&-\alpha_d&-\alpha_d&\cdots&-\alpha_d&1\\
1&-1&-1&-1&\cdots&-1&-1\end{array}\right)
\left(\begin{array}{c}
x_{0}\\ x_{1} \\ x_{2} \\x_{3}\\ \vdots \\x_{d-2}\\x_{d-1}\\x_{d}
\end{array}\right) = 
\left(\begin{array}{c}
 \alpha_{1} \\ \alpha_{2} \\\alpha_{3}\\ \alpha_{4} \\ \vdots \\ \alpha_{d-1}\\ \alpha_{d} \\1
\end{array}\right)
\end{equation*}}

Let $M_d$ be the $(d+1)\times(d+1)$ matrix of the system of equations 
and let $N_d$ be the
matrix that results from replacing the first column of $M_d$ with
the right-hand side of the system of equations. 
Then
$C_{112}(x,y,z)=\frac{\det(N_d)}{\det(M_d)}$. We start by
computing $\det(M_d)$. If we subtract $\alpha_j$ times the last
row from row $j$ for $j=1,2,\ldots,d$, (the elementary operation
is $R_j-\alpha_jR_{d+1}\rightarrow R_j$), and then subtract the
$j$th column from the $(j+1)$st column ($C_{j+1}-C_j\rightarrow
C_{j+1}$) for $j=2,3,\ldots,d$, and denote the resulting matrix by $M_{d}'$, then we get that
$\det(M_d)=\det(M'_d)$. 
Let $A_d$ be the matrix $M'_d$ without the first column and the
last row, and let $B_d$ be the matrix $M'_d$ without the second
column and the last row. Thus,
$\det(M_d)=(-1)^d(\det(A_d)+\det(B_d))$. It is easy to see that
$\det(A_d)=\prod_{j=1}^d(1+\alpha_j)$, and (by expanding along the
last row of $B_d$) that
$\det(B_d)=(1+\alpha_d)\det(B_{d-1})-(\alpha_d+\beta_d)\prod_{j=1}^{d-1}(1+\alpha_jy)$.
Using induction on $d$ with $\det(B_1)=-(\alpha_1+\beta_1)$ we get
that
$$\det(B_d)=-\sum_{j=1}^d(\alpha_j+\beta_j) \prod_{i=j+1}^d
(1+\alpha_i)\prod_{i=1}^{j-1}(1+\alpha_iy)=
-\prod_{i=1}^d(1+\alpha_i)\sum_{j=1}^d
{\textstyle \frac{\alpha_j+\beta_j}{1+\alpha_j}}\prod_{i=1}^{j-1}{\textstyle \frac{1+\alpha_iy}{1+\alpha_i}}.$$
Hence,
$$\det(M_d)=(-1)^d\prod_{i=1}^d(1+\alpha_i)\left(
1-\sum_{j=1}^d
\frac{\alpha_j+\beta_j}{1+\alpha_j}\prod_{i=1}^{j-1}\frac{1+\alpha_iy}{1+\alpha_i}\right).$$
Now, we consider the matrix $N_d$. If we add the first column to
all other columns, then it is easy to see that
$\det(N_d)=(-1)^d\prod_{i=1}^d(1+\alpha_i)$. Therefore, with  $\alpha_j=\frac{x^{2a_j}z^2}{1-x^{2a_j}z^2}$ and
$\beta_j=\frac{x^{a_j}z}{1+x^{a_j}z}$, we get (after algebraic simplification) that
$$\frac{\det(N_d)}{\det(M_d)}=\frac{1}{1-\sum_{j=1}^d
\left( x^{a_j}z\prod_{i=1}^{j-1}(1-(1-y)x^{2a_i}z^2)\right)}.$$
The case $|A|=\infty$ follows by taking the limit as $d \rightarrow \infty$. The proof for $C_{221}(x,y,z)$ follows with slight modifications.
\sof

Setting $y=0$ and $z = 1$ in Theorem~\ref{th112}, we
get that the generating function for the number of compositions of
$n$ with parts in $\NN$ that avoid the pattern $112$ is given by
$$C_{112}(x,0,1)= \frac{1}{1-\sum_{j\geq1}x^{j}\prod_{i=1}^{j-1}(1-x^{2i})}, $$
and the values of the corresponding sequence are 1, 1, 2, 4, 7,
13, 24, 43, 78, 142, 256, 463, 838, 1513, 2735, 4944, 8931, 16139,
29164, 52693, 95213 for $n=0,1,\ldots,20 $. 
The generating
function for the number of compositions of $n$ with parts in $\NN$
that avoid the pattern $221$ is given by
$$C_{221}(x,0,1)=\frac{1}{1-\sum_{i\geq1}\left(x^i\prod_{j\geq i+1}(1-x^{2j})\right)},$$
and the values of the corresponding sequence are 1, 1, 2, 4, 8, 15,
30, 58, 113, 220, 429, 835, 1627, 3169, 6172, 12023, 23419, 45616,
88853, 173073, 337118 for $n=0,1,\ldots,20$. Note that there are a lot less compositions of $n$ that avoid 112 than compositions that avoid $221$. This notion can be made more precise using the formulas for  the asymptotic behavior given in Theorem~\ref{asym}.

\subsection{The pattern $123$ (the statistic rise+rise)}
In the following theorem we will present the generating function
for the number of compositions of $n$ with $j$ parts in $A$ that
contain the pattern $123$ exactly $r$ times.

\begin{theorem}\label{th123}
Let $A$ be any ordered subset of $\NN$, with $|A|=d$. Then
$$C_{123}(x,y,z)=\frac{1}{1-t^1(A)-\sum\limits_{p=3}^d\sum\limits_{j=0}^{p-3}\binom{p-3}{j}t^{p+j}(A)(y-1)^{p-2}},$$
where $t^p(A)=\sum_{1\leq i_1<i_2<\cdots<i_p\leq d
}z^p\prod_{j=1}^p x^{a_{i_j}}$.
\end{theorem}
{\it Proof.}
Let $\sigma$ be any composition of $n$ with $m$ parts in
$A=\{a_1,\ldots,a_d\}$ that contains the pattern $123$ exactly
$r$ times. To derive recursions for the generating function, we
will break the composition into pieces, some of which have parts
in the set $A$, and others that have parts from the set
$A_k=\{a_{k+1},a_{k+2},\ldots,a_d\}=
A\backslash\{a_1,\ldots,a_k\}$ (the index of $A$ indicates the largest element excluded). To make this distinction for the
generating functions, we will indicate the specific set from which
the parts are selected as a superscript. 
Furthermore, since we want to split
off the parts $a_1, a_2,\ldots$ successively to create recursive
equations, we define $D^{A_k}(x,y,z)$ to be the generating
function for the number of compositions $\sigma$ of $n$ with $m$
parts in $A_k$ such that for $a \not\in A_k$, $a\sigma$ contains
the pattern $123$ exactly $r$ times.

For any composition $\sigma$ with parts in $A$, there are two
possibilities: either $\sigma$ does not contain $a_1$, in which
case the generating function is given by $C_{123}^{A_1}(x,y,z)$, or
the composition contains at least one occurrence of $a_1$, i.e., $\sigma
=\bar{\sigma}a_1\sigma_{k+1}\ldots\sigma_m$, where
${\bar{\sigma}}$ is a composition with parts from  ${A_1}$, with
generating function $C_{123}^{A_1}(x,y,z)C_{123}^A(a_1|x,y,z)$.
Altogether, we have
\begin{equation}\label{C123}
C_{123}^A(x,y,z)=C_{123}^{A_1}(x,y,z) +C_{123}^{A_1}(x,y,z)C_{123}^A(a_1|x,y,z).
\end{equation}
Now let us consider the compositions $\sigma$ of $n$ with $m$
parts in $A$ starting with $a_1$ which contain the pattern $123$
exactly $r$ times. Again, there are two cases: either $\sigma$
contains exactly one occurrence of $a_1$, or the part $a_1$ occurs
at least twice in $\sigma$. In the first case, the generating
function is given by $x^{a_1}z\,D^{A_1}(x,y,z)$. 
If $\sigma$ contains $a_1$ at least twice, then we split
the composition into pieces according to the second occurrence of
$a_1$, i.e., $\sigma 
=a_1\bar{\sigma}a_1\sigma_{k+1}\ldots\sigma_m$, where
${\bar{\sigma}}$ is a (possibly empty) composition with parts from
$A_1$. Splitting off the initial part $a_1$ results in the
generating function
$x^{a_1}z\,D^{A_1}(x,y,z)\,C_{123}^{A}(a_1|x,y,z)$.  Thus,
$$C_{123}^A(a_1|x,y,z)=x^{a_1}zD^{A_1}(x,y,z)+x^{a_1}zD^{A_1}(x,y,z)C_{123}^A(a_1|x,y,z).$$
Solving for $C_{123}^A(a_1|x,y,z)$ and 
substituting  into (\ref{C123}) gives
\begin{equation}\label{CD123}
C_{123}^A(x,y,z)=\frac{C_{123}^{A_1}(x,y,z)}{1-x^{a_1}zD^{A_1}(x,y,z)}.
\end{equation}
We now derive an expression for $D^{A_1}(x,y,z)$ by considering
compositions $\sigma$ with parts in $A_1$ such that $a_1\sigma$
contains the pattern 123 exactly $r$ times. If $\sigma$ does not contain the part $a_2$, the
generating function for $\sigma$ is given by $D^{A_2}(x,y,z)$. Otherwise, we  write $\sigma =
\bar{\sigma}^1a_2\bar{\sigma}^2a_2\bar{\sigma}^3\ldots
a_2\bar{\sigma}^{\ell+2}$ with $\ell \ge 0$, where
$\bar{\sigma}^{j}$ is a (possibly empty) composition with parts in
$A_2$ for $j=1,\ldots,\ell+2$. There are four subcases, depending
on whether $\bar{\sigma}^1$ and $\bar{\sigma}^2$ are empty
compositions or not. If $\bar{\sigma}^1$ and $\bar{\sigma}^2$ are
both empty compositions, then $a_1\sigma=a_1a_2$
or $a_1\sigma=a_1a_2a_2\bar{\sigma}^3\ldots a_2\bar{\sigma}^{\ell+2}$, $\ell
\ge 1$. In either case we can split off the initial part $a_2$ of $\sigma$
which results in one less part, but no reduction in the
occurrences of the pattern 123. Thus, the generating function for
$\sigma$ is given by $$x^{a_2}z\sum_{\ell\geq0}
(x^{a_2}z\,D^{A_2}(x,y,z))^{\ell}=\frac{x^{a_2}z}{1-x^{a_2}z\,D^{A_2}(x,y,z)}.$$
If $\bar{\sigma}^1$ is the empty composition and $\bar{\sigma}^2$
is not the empty composition, then $a_1\sigma =
a_1a_2\bar{\sigma}^2$ or $a_1\sigma =
a_1a_2\bar{\sigma}^2a_2\bar{\sigma}^3\ldots
a_2\bar{\sigma}^{\ell+2}$, $\ell \ge 1$, and  the generating function for
$\sigma$ is given by
$$x^{a_2}z\,y(D^{A_2}(x,y,z)-1)\sum_{\ell\geq0}(x^{a_2}zD^{A_2}(x,y,z))^{\ell}
=\frac{x^{a_2}z\,y(D^{A_2}(x,y,z)-1)}{1-x^{a_2}zD^{A_2}(x,y,z)}.$$

If $\bar{\sigma}^1$ is not the empty composition and
$\bar{\sigma}^2$ is the empty composition, then $\sigma =
\bar{\sigma}^1a_2$ or $\sigma =
\bar{\sigma}^1a_2a_2\bar{\sigma}^3\ldots
a_2\bar{\sigma}^{\ell+2}$, $\ell \ge 1$, and
the generating function for $\sigma$  is given by
$$x^{a_2}z(D^{A_2}(x,y,z)-1)\sum_{\ell\geq0}
(x^{a_2}zD^{A_2}(x,y,z))^{\ell}=\frac{x^{a_2}z(D^{A_2}(x,y,z)-1)}
{1-x^{a_2}zD^{A_2}(x,y,z)}.$$
Finally, if both $\bar{\sigma}^1$ and $\bar{\sigma}^2$ are nonempty
compositions, then the generating function for $\sigma$ is given by
$$\frac{x^{a_2}z(D^{A_2}(x,y,z)-1)^2}{1-x^{a_2}zD^{A_2}(x,y,z)}.$$

Adding all four cases, using the shorthand $D^A$ for  $D^A(x,y,z)$, and solving for $D^{A_{1}}$ we get
\begin{equation}\label{DA1}
D^{A_1}=\frac{(1-x^{a_2}z(1-y))D^{A_2}+x^{a_2}z(1-y)}{1-x^{a_2}zD^{A_2}}.
\end{equation}

We give an explicit expression for $D^A$ in the following lemma.

\begin{lemma} \label{lem1} Let $A=\{a_1,\ldots,a_d\}$ and
$t^p(A_k)=\sum_{k+1\leq i_1<i_2<\cdots<i_p\leq d
}z^p\prod_{j=1}^p x^{a_{i_j}}$ for all $p$ and $k=0,1,\ldots,d$.
Then
\begin{equation}\label{DA}
D^A=\frac{1+\sum_{p=2}^d\sum_{j=0}^{p-2}\binom{p-2}{j}t^{p+j}(A)
(y-1)^{p-1}}{1-t^1(A)-\sum_{p=3}^d\sum_{j=0}^{p-3}\binom{p-3}{j}t^{p+j}(A)(y-1)^{p-2}}.
\end{equation}
\end{lemma}
{\it Proof.}
Before we start proving (\ref{DA}), we will give an interpretation of $t^p(A_k)$ as the generating function for the
number of partitions with $p$ distinct parts from the set $A_k$,
where $A_k$ has $d-k$ elements, with $A_0=A$ and $0 \le k \le d-1$. These partitions either contain the
part $a_{k+1}$ or not. In the first case, the generating function is
given by $x^{a_{k+1}}z\,t^{p-1}(A_{k+1})$, and in the second case, by
$t^p(A_{k+1})$. Thus,
\begin{equation} \label{tpd}
t^p(A_k)=t^p(A_{k+1})+x^{a_{k+1}}z\,t^{p-1}(A_{k+1}).
\end{equation}

We now prove (\ref{DA}) by induction on $d$, the number of elements in
$A$. For $d=0$ and $d=1$ we have $D^{\emptyset}=1$ and
$D^{\{a_1\}}=\sum_{n\ge0}x^{n\,a_1}z^n=1/(1-x^{a_1}z)$,
respectively, and thus, (\ref{DA}) holds. Now assume that (\ref{DA}) holds for $d-1$. Using (\ref{DA1}) and
the induction hypothesis for the set $A_1$ gives
$$\begin{array}{ll}
D^A&=\frac{(1-x^{a_1}z(1-y))D^{A_1}+x^{a_1}z(1-y)}{1-x^{a_1}zD^{A_1}}\\
&{\small
=\frac{(1-x^{a_1}z(1-y))\left(1+\sum\limits_{p=2}^{d-1}\sum\limits_{j=0}^{p-2}\binom{p-2}{j}
t^{p+j}(A_1)(y-1)^{p-1}\right)}
{1-t^1(A_1)-\sum\limits_{p=3}^{d-1}\sum\limits_{j=0}^{p-3}\binom{p-3}{j}t^{p+j}(A_1)(y-1)^{p-2}
-x^{a_1}z\left(1+\sum\limits_{p=2}^{d-1}\sum\limits_{j=0}^{p-2}\binom{p-2}{j}t^{p+j}(A_1)(y-1)^{p-1}\right)}}\\
&{\small\quad +\quad
\frac{x^{a_1}z(1-y)\left(1-t^1(A_1)-\sum\limits_{p=3}^{d-1}\sum\limits_{j=0}^{p-3}
\binom{p-3}{j}t^{p+j}(A_1)(y-1)^{p-2}\right)}
{1-t^1(A_1)-\sum\limits_{p=3}^{d-1}\sum\limits_{j=0}^{p-3}\binom{p-3}{j}t^{p+j}(A_1)(y-1)^{p-2}
-x^{a_1}z\left(1+\sum\limits_{p=2}^{d-1}\sum\limits_{j=0}^{p-2}\binom{p-2}{j}t^{p+j}(A_1)(y-1)^{p-1}\right)}}=\dfrac{s_1}{s_2}.
\end{array}$$
We first rewrite the denominator and obtain
\begin{align*}
s_{2}&=1-t^1(A_1)-x^{a_1}z-\sum\limits_{p=3}^{d-1}\sum\limits_{j=0}^{p-3}\tbinom{p-3}{j}t^{p+j}(A_1)(y-1)^{p-2}\\
&\phantom{====}-\sum\limits_{p=2}^{d-1}\sum\limits_{j=0}^{p-2}\tbinom{p-2}{j}x^{a_1}z\,t^{p+j}(A_1)(y-1)^{p-1}.
\end{align*}
Combining the two double sums by reindexing the second one, adding the terms for $p=d$ to the first one (which by definition are all zero, as $t_{d-1}^{d+j}(A_1)=0$ for $j \ge 0$) and applying (\ref{tpd}) gives
\vspace{-10pt}
$$s_2=1-t^1(A)
-\sum\limits_{p=3}^{d}\sum\limits_{j=0}^{p-3}\tbinom{p-3}{j}t^{p+j}(A)(y-1)^{p-2}.$$  Next we rewrite the numerator:
\begin{align*}
s_{1} &=1+x^{a_1}z(y-1)t^1(A_1)+\sum\limits_{p=2}^{d-1}\sum\limits_{j=0}^{p-2}\tbinom{p-2}{j}
t^{p+j}(A_1)(y-1)^{p-1}\\
&{}\phantom{=1}+ \,x^{a_1}z\sum\limits_{p=2}^{d-1}\sum\limits_{j=0}^{p-2}\tbinom{p-2}{j}
t^{p+j}(A_1)(y-1)^{p}+ \,x^{a_1}z\sum\limits_{p=3}^{d-1}\sum\limits_{j=0}^{p-3}
\tbinom{p-3}{j}t^{p+j}(A_1)(y-1)^{p-1}.
\end{align*}
We now look at the coefficient of $(y-1)^m$ and collect terms according to the respective power. 
For $m=0$, the coefficient is $1$. If $m=1$, then the coefficient is given by $x^{a_1}z\,t^1(A_1)+t^2(A_1)=t^2(A)$ (using (\ref{tpd})). If $m=2,3,\ldots,d-1$, then the coefficient of $(y-1)^m$ is equal to
\begin{align*}
{}&{}\sum_{j=0}^{m-1}\tbinom{m-1}{j}t^{m+1+j}(A_1)+x^{a_1}z\sum_{j=0}^{m-2}\tbinom{m-2}{j}t^{m+j}(A_1)
+x^{a_1}z\sum_{j=0}^{m-2}\tbinom{m-2}{j}t^{m+1+j}(A_1)\\
&=\sum_{j=0}^{m-1}\tbinom{m-1}{j}t^{m+1+j}(A_1)+x^{a_1}z\left(\sum_{j=0}^{m-2}\tbinom{m-2}{j}t^{m+j}(A_1)
+\sum_{j=1}^{m-1}\tbinom{m-2}{j-1}t^{m+j}(A_1)\right)\\
\intertext{which, using the identity $\binom{a}{b-1}+\binom{a}{b}=\binom{a+1}{b}$ and the fact that $\binom{m-2}{-1}=\binom{m-2}{m-1}=0$,}
&=\sum_{j=0}^{m-1}\tbinom{m-1}{j}t^{m+1+j}(A_1)+x^{a_1}z\sum_{j=0}^{m-1}\tbinom{m-1}{j}t^{m+j}(A_1)\\
&=\sum_{j=0}^{m-1}\tbinom{m-1}{j}\left(t^{m+1+j}(A_1)+x^{a_1}z\,t^{m+j}(A_1)\right)
=\sum_{j=0}^{m-1}\tbinom{m-1}{j}t^{m+j+1}(A),
\end{align*}
where the last equality follows once more from using (\ref{tpd}). Thus,  \begin{align*}
s_1 &= 1 + t^2(A)(y-1) + \sum_{m=2}^{d-1}\sum_{j=0}^{m-1}\tbinom{m-1}{j}t^{m+j+1}(A)(y-1)^m\\
&=1 + \sum_{m=1}^{d-1}\sum_{j=0}^{m-1}\tbinom{m-1}{j}t^{m+j+1}(A)(y-1)^m.
\end{align*}
Reindexing the sum and combining this result with the result for $s_2$ completes the proof of the lemma.\sof

We now can obtain an exact formula for
$C_{123}^A(x,y,z)$ as follows. Using (\ref{CD123}) and Lemma~\ref{lem1}
results in
$$C_{123}^A(x,y,z)=\frac{C_{123}^{A_1}(x,y,z)}{1-x^{a_1}z
\frac{1+\sum_{p=2}^{d-1}\sum_{j=0}^{p-2}\binom{p-2}{j}t^{p+j}(A_1)(y-1)^{p-1}}
{1-t^1(A_1)-\sum_{p=3}^{d-1}\sum_{j=0}^{p-3}\binom{p-3}{j}t^{p+j}(A_1)(y-1)^{p-2}}}.$$

Using the same arguments as in the proof of the above lemma
we get that
\begin{align*}\quad C_{123}^A(x,y,z)&=C_{123}^{A_1}(x,y,z)\frac{1-t^1(A_1)-\sum_{p=3}^{d-1}\sum_{j=0}^{p-3}\binom{p-3}{j}t^{p+j}(A_1)(y-1)^{p-2}}{
1-t^1(A)-\sum_{p=3}^{d}\sum_{j=0}^{p-3}\binom{p-3}{j}t^{p+j}(A)(y-1)^{p-2}}.\phantom{---}\\
\intertext{Iterating this equation and using that $C_{123}^{A_d}(x,y,z)=C_{123}^{\emptyset}(x,y,z)=1$ results in:}
\quad C_{123}^A(x,y,z)&=\prod_{k=0}^{d-1}\frac{1-t^1(A_{k+1})-\sum_{p=3}^{d-1}\sum_{j=0}^{p-3}\binom{p-3}{j}t^{p+j}(A_{k+1})(y-1)^{p-2}}{
1-t^1(A_k)-\sum_{p=3}^{d}\sum_{j=0}^{p-3}\binom{p-3}{j}t^{p+j}(A_k)(y-1)^{p-2}}.\\
\end{align*} Simplifying and using that
$t^p(A_d)=t^p(\emptyset)=0$ gives the desired result.
\sof

To apply Theorem~\ref{th123} to $A=\NN$, we first show that $t^p(\NN)=x^{\binom{p+1}{2}}z^p/(x;x)_{p}$, where we use the customary notation $(a;q)_{n}=\prod_{{j=0}}^{n-1}(1-a\,q^{j})$. Clearly, the formula for $t^p(\NN)$ holds for $p=0$. For  $p\geq1$, we get from the definition of $t^p(\NN)$ with $a_i = i$ for $i \ge 1$ that
\begin{align*}
t^p(\NN) &= z^p\sum_{1\leq
i_1<i_2<\cdots<i_{p-1}<i_p}x^{i_1+i_2+\cdots+i_{p}}= z^p\sum_{1\leq
i_1<i_2<\cdots<i_{p-1}}x^{i_1+i_2+\cdots+i_{p-1}}\sum_{i >i_{p-1}}x^i\\
&= z^p\sum_{1\leq
i_1<i_2<\cdots<i_{p-1}}x^{i_1+i_2+\cdots+i_{p-1}}\frac{x^{i_{p-1}+1}}{(1-x)}\\
&= z^p\frac{x}{(1-x)}\sum_{1\leq i_1<i_2<\cdots<i_{p-2}}x^{i_1+i_2+\cdots+i_{p-2}}\sum_{i>i_{p-2}}x^{2i}\\
&= \dots=x^{\binom{p+1}{2}}z^p\prod_{j=1}^p (1-x^j)^{-1}=x^{\binom{p+1}{2}}z^p/(x;x)_{p}.
\end{align*}

Setting $y=0$ and $z=1$ in Theorem~\ref{th123}  we obtain the generating function for the number of compositions with parts in $\NN$ that avoid 123 as
$$C_{123}^\NN(x,0,1)=\frac{1}{1-\frac{x}{1-x}-\sum\limits_{p\geq3}\sum\limits_{j=0}^{p-3}
\binom{p-3}{j}\frac{x^{\binom{p+1+j}{2}}}{(x;x)_{p+j}}(-1)^{p-2}},$$
and the sequence for the number of 123 pattern-avoiding compositions with parts in $\NN$ for $n=0$ to $n=20$ is given by $1$, $1$, $2$, $4$, $8$, $16$, $31$, $61$, $119$, $232$, $453$, $883$, $1721$, $3354$, $6536$, $12735$, $24813$, $48344$, $94189$, $183506$, $357518$. Note that the first time the pattern 123 can occur is for $n=6$, as the composition 123.

\subsection{The patterns $\{121,132,231\}$ and $\{212, 213, 312\}$ (the statistics $peak$ = rise + drop and  $valley$ = drop+rise)}\label{sec peak}

We will now look at the set of patterns $\{121,132,231\}$ together,
as they constitute the statistic $peak$, and likewise for the set  $\{212, 213, 312\}$. For ease of use, we will
refer to these sets of patterns collectively as the patterns $peak$ and $valley$, respectively,
and define $C_{peak}^A(x,y,z)$ and $C_{valley}^A(x,y,z)$ accordingly.  Before we can state the result for the respective generating functions, we need a few definitions. For any set $B \subseteq A$ 
and for $s\geq1$, we define
$$P^{s}(B)=\{(i_1,\ldots,i_{s}) |\, a_{i_{j}}\in B, j=1,\ldots,s, \mbox{ and } i_{2\ell-1} < i_{2\ell} \le i_{2\ell+1} \mbox{ for } 1 \le \ell \le \lfloor s/2\rfloor\},$$
$$Q^{s}(B)=\{(i_1,\ldots ,i_{s}) |\, a_{i_{j}}\in B, j=1,\ldots,s, \mbox{ and } i_{2\ell-1} \le i_{2\ell} < i_{2\ell+1} \mbox{ for } 1 \le \ell \le \lfloor s/2\rfloor \}$$
and 
$$M^s(B)=\sum_{(i_1,\ldots,i_s)\in
P^s(B)} z^{p} \prod_{j=1}^s b_{i_j} \quad \mbox{and } \quad N^s(B)=\sum_{(i_1,\ldots,i_s)\in Q^s(B)} z^{p} \prod_{j=1}^s
b_{i_j}.  $$

On route to the explicit expressions given in Theorem~\ref{peak} we will express the generating functions for the patterns $peak$ and $valley$ as continued fractions, for which we will use the notation { $[c_0,c_1,c_2,\ldots,c_{n-1},c_n]=c_0+\cfrac{1}{c_1+\cfrac{1}{c_2+\cfrac{1}{\ddots+\cfrac{1}{c_{n-1}+1/c_n}}}}.$}

\begin{theorem}\label{peak}
Let $A=\{a_1,\ldots, a_d \}$,  $P^{s}(A)$, $Q^{s}(A)$, $M^s(A)$, and $N^s(A)$ defined as
above. Then 
$$C_{peak}^A(x,y,z)=\frac{1+\sum_{j\ge 1}M^{2j}(A) (1-y)^j}{1+\sum_{j\ge 1}M^{2j}(A) (1-y)^j-\sum_{j\ge 0}M^{2j+1}(A) (1-y)^{j}}, \quad \mbox{and}$$
$$C_{valley}^A(x,y,z)=\frac{1+\sum_{j\ge 1}M^{2j}(A) (1-y)^{j}}{1+\sum_{j\ge 1}M^{2j}(A)
(1-y)^{j}-\sum_{j\ge 0}N^{2j+1}(A) (1-y)^{j}}.$$
\end{theorem}
{\it Proof.}
We prove the result for the pattern $peak$.  To derive first a recursion, and then an explicit formula for $C_{peak}^A(x,y,z)$,
we concentrate on occurrences of $a_d$, the largest part in the set $A=\{a_1, a_2, \ldots,a_d\}$.
If $a_d$ is surrounded by smaller parts on both sides, then a peak occurs. Let $\sigma$ be
any composition with parts in $A$, and define
$\bar{A}_k= \{a_1,\ldots,a_k\}$
(the index for the set indicates the largest element included). Note that $A = \bar{A}_d$. \\
We now look at the different possibilities for occurrences of
$a_d$. If $\sigma$ does not contain $a_d$, then the generating
function is given by $C_{peak}^{\bar{A}_{d-1}}(x,y,z)$. If there
is at least one occurrence, then we need to look at three cases
for the first occurrence of $a_d$. If the first occurrence is at
the beginning of the composition, i.e., $\sigma = a_d\sigma'$,
where $\sigma'$ is a (possibly empty) composition with parts in $A$, then no peak occurs, and the generating function is given
by $x^{a_d}z\,C_{peak}^{A}(x,y,z)$. If the first (and only)
occurrence of $a_d$ is at the end of the composition, i.e.,
$\sigma = \bar{\sigma}a_d$, where $\bar{\sigma}$ is a non-empty
composition with parts in $\bar{A}_{d-1}$, then the generating
function is given by
$x^{a_d}z(C_{peak}^{\bar{A}_{d-1}}(x,y,z)-1)$. Finally, if the
first occurrence of $a_d$ is in the interior of the composition,
then $\sigma = \bar{\sigma}a_d\sigma'$, where $\sigma'$ is a
non-empty composition with parts in $A$. If $\sigma'$ starts with
$a_d$, then $\sigma = \bar{\sigma}a_d\,a_d\,\sigma'$, and both
$a_d$'s can be split off without decreasing the number of
occurrences of peaks; the generating function is given by
$(C_{peak}^{\bar{A}_{d-1}}(x,y,z)-1)x^{2a_d}z^2C_{peak}^{A}(x,y,z)$.
If $\sigma'$ does not start with $a_d$, then a peak occurs and the
generating function is given by
$$ (\underbrace{C_{peak}^{\bar{A}_{d-1}}(x,y,z)-1}_{\bar{\sigma}\mbox{ \Small non-empty}})x^{a_d}z\,y(\underbrace{C_{peak}^{A}(x,y,z)-1}_{\sigma'\mbox{ \Small non-empty}}\underbrace{-x^{a_d}z\,C_{peak}^{A}(x,y,z)}_{\mbox{\Small does not start with }a_d}).$$
Now, let $C^A=C_{peak}^A(x,y,z)$. Combining the three cases above, we get
\begin{align*}
C^A&=C^{\bar{A}_{d-1}}+x^{a_d}z\,C^A+x^{a_d}z\,(C^{\bar{A}_{d-1}}-1)+x^{2a_d}z^2C^A(C^{\bar{A}_{d-1}}-1)\\
&\phantom{==}+{}x^{a_d}zy(C^A-1-x^{a_d}zC^A)(C^{\bar{A}_{d-1}}-1),\end{align*}
or, equivalently,
\begin{align}\label{Cafrac}
C^A&=\dfrac{(1+x^{a_d}z(1-y))C^{\bar{A}_{d-1}}-x^{a_d}z(1-y)}{1-x^{a_d}z(1-x^{a_d}z)(1-y)-x^{a_d}z(x^{a_d}z(1-y)+y)C^{\bar{A}_{d-1}}}\nonumber \\
&=\dfrac{1}{1-x^{a_d}z-\dfrac{C^{\bar{A}_{d-1}}-1}{(1+x^{a_d}z(1-y))C^{\bar{A}_{d-1}}-x^{a_d}z(1-y)}} \\
&=\cfrac{1}{1-x^{a_d}z-\cfrac{1}{[x^{a_d}z(1-y),1-1/C^{\bar{A}_{d-1}}]}}. \nonumber 
\end{align}
Hence, by induction on $d$ and using the fact that
$C^{\bar{A}_1}=\frac{1}{1-x^{a_1}z}$, we can express the
generating function $C_{peak}^A(x,y,z)$ as a continued fraction. 

\begin{lemma}\label{cf}
For $A=\{a_1,\ldots,a_d\}$, $b_i=x^{a_i}z$, and $C^A= C^{A}_{peak}(x,y,z)$,
$$C^{A}=
\cfrac{1}{1-b_d-\cfrac{ 1}{[b_d(1-y),b_{d-1},b_{d-1}(1-y),\ldots,b_2,b_2(1-y),b_1]}}.$$
\end{lemma}
Now we derive an explicit formula for $C^A_{peak}(x,y,z)$ based on
recursions for $M^s(\bar{A}_{d})=M^s(A)$ for odd and even $s$; if $s$ is odd, both
the last and second-to-last elements can equal $d$, whereas in the
case $s$ even, the second-to-last element can be at most $d-1$. By
separating the elements of $P^s(A)$ according to whether the last
element equals $d$ or is less than $d$, we get the following two
recursions:
\begin{align} \label{md}
M^{2s+1}(A) &=b_d\,M^{2s}(A)+M^{2s+1}(\bar{A}_{d-1}) \quad \mbox{and} \\
M^{2s}(A) &=b_d\, M^{2s-1}(\bar{A}_{d-1})+M^{2s}(\bar{A}_{d-1}).\nonumber
\end{align}
Define
$G_d=
\frac{ 1}{[b_d(1-y),b_{d-1},b_{d-1}(1-y),\ldots,b_2,b_2(1-y),b_1]} $,  i.e., $G_d$ consists of the portion of $C^A$ in the continued fraction expansion which has a repeating pattern. We now derive an expression for $G_d$ in terms of the $M^s(A)$.

\begin{lemma} \label{gfrac}
For all $d\geq2$, $G_d=\dfrac{\sum_{j\ge 0} M^{2j+1}(A)
(1-y)^{j}}{1+\sum_{j\ge 1}M^{2j}(A) (1-y)^j}=\dfrac{G_d^1}{G_d^2}$.
\end{lemma}
{\it Proof.}
We prove the statement by induction on $d$. For $d=2$, $G_2^1=b_1$
(only the term $j = 0$ contributes), $G_2^2=1+b_1\,b_2\,(1-y)$
(only the term $j = 1$ contributes) and therefore, the lemma
holds. Now let $d\geq3$ and assume that the lemma holds for $d-1$,
i.e., $G_{d-1}=\frac{G_{d-1}^1}{G_{d-1}^2}$. By the definition of
$G_d$ it is easy to see that $G_{d}=\frac{1}{b_{d}(1-y)+1/(b_{d-1}+G_{d-1})}.$
Substituting the induction hypothesis for $d-1$ into the expression for $G_d$ and simplifying yields
$$G_{d}=\frac{b_{d-1}G_{d-1}^2+G_{d-1}^1}{G_{d-1}^2+b_d(1-y)(b_{d-1}G_{d-1}^2+G_{d-1}^1)}.$$
Using the definitions of $G_{d-1}^1$ and $G_{d-1}^2$ and (\ref{md}) yields $b_{d-1}G_{d-1}^2+G_{d-1}^1=G_d^1$, and this result together with the definitions of $G_{d-1}^1$ and $G_{d-1}^2$ and (\ref{md}) yields $G_{d-1}^2+b_d(1-y)(b_{d-1}G_{d-1}^2+G_{d-1}^1)=G_d^2$, which completes the proof of Lemma~\ref{gfrac}. \sof

Now we use Lemma~\ref{gfrac} and $C^A=\frac{1}{1-b_d-G_d}$ to get (after simplification) that $C^{A}$ equals 
\begin{align*}
\frac{1+\sum_{j\ge 1}M^{2j}(A) (1-y)^j}{1-b_d-M^1(A)+\sum_{j\ge 1}M^{2j}(A) (1-y)^j-\sum_{j\ge 1}(b_dM^{2j}(A)+M^{2j+1}(A)) (1-y)^{j}}.
\end{align*}
Using (\ref{md}) we get that
\begin{align*}
C^A&=\frac{1+\sum_{j\ge 1}M^{2j}(A) (1-y)^j}{1-M^1(A)+\sum_{j\ge 1}M^{2j}(A) (1-y)^j-\sum_{j\ge 1}M^{2j+1}(A) (1-y)^{j}}\\
&=\frac{1+\sum_{j\ge 1}M^{2j}(A) (1-y)^j}{1+\sum_{j\ge 1}M^{2j}(A) (1-y)^j-\sum_{j\ge 0}M^{2j+1}(A) (1-y)^{j}}.
\end{align*}
This completes the proof for the pattern $peak$ since the formula also holds for the case $d=\infty$. The result for the pattern $valley$ follows with minor modifications, focusing on the smallest rather than the largest part, replacing $\bar{A}_{k}$ with $A_{k}=\{a_{k+1},a_{k+2},\ldots,a_{d}\}$, and using the  recursions 
$$M^s(A_k)=b_{k+1}\,N^{s-1}(A_{k+1})+M^{s}(A_{k+1}) \mbox{  and  }
N^s(A_k) =b_{k+1}M^{s-1}(A_k)+N^s(A_{k+1}), $$ which are obtained by
separating the elements of $P^s(A_k)$ according to whether
the first element equals $k+1$ or is greater than $k+1$.
\sof

We now apply Theorem~\ref{peak} to $A=\NN$. Similar to the derivation of $t^{p}(\NN)$, but with extra care since there are both $<$ and $\le$ constraints (resulting in only odd powers of $x$ in the numerator), we get for all $s\geq1$,
\begin{align*}
M^{2s}(\NN)&=\sum_{1\leq i_1<i_2\leq i_{3}\cdots i_{2s-2}\leq i_{2s-1}<i_{2s}}x^{i_1+\cdots+i_{2s}}z^{2s}
=x^{s(s+2)}z^{2s}/(x;x)_{{2s}},\\
M^{2s+1}(\NN)&=\sum_{1\leq i_1<i_2\leq i_{3}\cdots i_{2s-1}<i_{2s}\leq i_{{2s+1}}}x^{i_1+\cdots+i_{2s}}z^{2s}
=x^{s^2+3s+1}z^{2s+1}/(x;x)_{{2s+1}}, \mbox{ and}\\
N^{2s+1}(\NN)&=\sum_{1\leq i_1\leq i_2<i_{3}\cdots i_{2s-1}\leq i_{2s}< i_{{2s+1}}}x^{i_1+\cdots+i_{2s+1}}z^{2s+1}=x^{(s+1)^2}z^{2s+1}/(x;x)_{{2s+1}}.
\end{align*}
Substituting these expressions into Theorem~\ref{peak} gives
the generating functions for the number of compositions of $n$
with parts in $\NN$ without $peaks$ and $valleys$, respectively, as
$$C_{peak}^\NN(x,0,1)=\frac{1+\sum_{j\ge 1}\frac{x^{j(j+2)}}{(x;x)_{{2j}}}}{1+\sum_{j\ge 1}\frac{x^{j(j+2)}}{(x;x)_{{2j}}}
-\sum_{j\ge0}\frac{x^{j^2+3j+1}}{(x;x)_{{2j+1}}}}$$
and 
$$C_{valley}^\NN(x,0,1)=\frac{1+\sum_{j\ge 1}\frac{x^{j(j+2)}}{(x;x)_{{2j}}}}{1+\sum_{j\ge 1}\frac{x^{j(j+2)}}{(x;x)_{{2j}}}
-\sum_{j\ge0}\frac{x^{(j+1)^2}}{(x;x)_{{2j+1}}}}.$$
Rewriting the generating function as a geometric series allows us to compute the sequence for the number of $peak$-avoiding
compositions with parts in $\NN$. The terms for $n=0$ to $n=20$
are given by $1$, $1$, $2$, $4$, $7$, $13$, $22$, $38$, $64$,
$107$, $177$, $293$, $481$, $789$, $1291$, $2110$, $3445$, $5621$,
$9167$, $14947$, $24366$. The corresponding sequence for $valley$-avoiding compositions is given by $1$, $1$, $2$, $4$, $8$, $15$, $28$, $52$, $96$,
$177$, $326$, $600$, $1104$, $2032$, $3740$, $6884$, $12672$,
$23327$, $42942$, $79052$, $145528$. Note that the first time a
peak can occur is for $n=4$ (121), and the first time a
valley can occur is for $n=5$ (212).

{\bf Remark:} Even though the statistics $peak$ and $valley$
are in some sense symmetric, one cannot obtain the number of
$valley$-avoiding compositions from the number of $peak$-avoiding
compositions. However, there is a connection, namely the number of
$valleys$ in the compositions of $n$ with $m$ parts are equal to the
number of $peaks$ in the compositions of $m(n+1)-n$ with $m$ parts.
This can easily be seen as follows: In each composition of $n$
with $m$ parts, replace each part $\sigma_i$ by $(n+1)-\sigma_i$,
which results in a composition of
$m(n+1)-\sum_{i=1}^m\sigma_i=m(n+1)-n$. This connection will be important in Section~\ref{words}, when we apply the results derived for the various patterns to words on $k$ letters.

\vskip 30pt

\section{Asymptotics for the number of compositions avoiding $\tau$} \label{asymps}

We will now use methods from Complex Analysis to compute the asymptotics for the number of compositions with parts in $\NN$ which avoid a given pattern $\tau$. We think of the generating function as a complex function, and indicate this fact by using the variable $z$ instead of the variable $x$. Thus, we look at the function $C_{\tau}(z)=C_{\tau}(z,0,1)=\sum_{n \ge 0}C_{\tau}(n,0)z^n$. Since $C_{\tau}(z)$  is meromorphic,  the asymptotic behavior of $C_{\tau}(n,0)$ is determined by the dominant pole of the function $C_{\tau}(z) = 1/f(z)$, i.e., the smallest positive root of $f(z)$ (see for example~\cite{Wilf1990}). Using Theorem 5.2.1~\cite{Wilf1990} and the discussions preceding it, we obtain the following result.

\begin{theorem} \label{asym} The asymptotic behavior for $\tau$-avoiding compositions with parts in $\NN$ is given by \begin{eqnarray*}
  C_{111}(n,0) &= & 0.499301 \cdot 1.91076^n+O((10/7)^n) \\
  C_{112}(n,0)&= & 0.692005 \cdot 1.80688^n+O((10/7)^n)\\
  C_{221}(n,0)&= & 0.545362 \cdot 1.94785^n+O((10/7)^n)\\
  C_{123}(n,0) &= & 0.576096 \cdot 1.94823^n+O((10/7)^n)\\
  C_{peak}(n,0)&= & 1.394560 \cdot 1.62975^n+O((10/7)^n)\\
  C_{valley}(n,0)&= &  0.728207 \cdot 1.84092^n+O((10/7)^n).
 \end{eqnarray*}
 \end{theorem}
{\it Proof.} Let $\rho$ be the smallest positive root of $f(z)$. If $\rho$ is a simple pole, then the residue is given by $1/ f'(\rho)$. Since $C_{\tau}(n,0) \le 2^{n-1}$, the number of unrestricted compositions with parts in $\NN$, we know that the radius of convergence of $C_{\tau}(z)>0.5$, and therefore, $\rho > 0.5$ for all patterns $\tau$. Using Mathematica and Maple, we compute both $\rho$ and $1/ f'(\rho)$ for all patterns. To verify that we are dealing with simple poles in each case, we use the {``Principle of the Argument''} (see Theorem 4.10a,  \cite{Hen1974}), which states that the number of zeros of a function $f(z)$ is equal to the winding number of the
transformed curve $f(\Gamma)$ around the origin, where  $\Gamma$ is a simple closed curve. We use as $\Gamma$ the circle $r=|z|=0.7$. Figure~1 shows the six graphs. 
\begin{center}
\begin{figure}[ht]\label{allgraphs}
\vspace{-40pt} 
\epsfxsize=350.0pt \epsffile{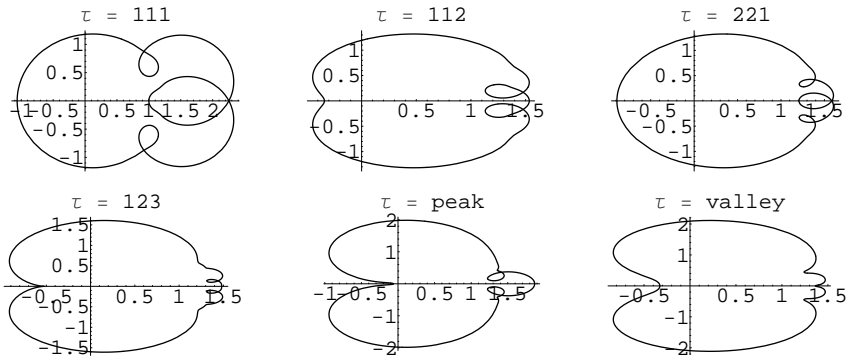}
\caption{The image of the circle $|z|=0.7$ under the respective generating functions.}
\end{figure}
\end{center}
\vspace{-40pt} Clearly, the winding number is 1 in each case, confirming that $\rho$ is a simple pole. Thus, we obtain from Theorem 5.2.1~\cite{Wilf1990} that  $C_{\tau}(n,0)=K\cdot v^n +O((1/r)^n)$, where $K=-1/(\rho f'(\rho))$ and $v=1/\rho$, which completes the proof. \sof

Note that Theorem~\ref{asym} for $\tau = 111$ gives the asymptotics for the number of $2$-Carlitz compositions. Asymptotics for the Carlitz compositions were given in~\cite{KnoPro1998}.

\vskip 30pt


\section{Counting occurrences of subword patterns} \label{words}

Several authors (for example, see~\cite{BurMan2003, BurMan2003a} and references
therein) have studied the occurrence of subword patterns in words
on $k$ letters. We will apply the results derived in the previous
sections to words on $k$ letters, and thus obtain previous results
as special cases.

Let $[k]=\{1,2,\dots,k\}$ be a (totally ordered) alphabet on $k$
letters. We call the elements of $[k]^n$ \emph{words}. A word
$\sigma$ \emph{contains} a pattern $\tau$ if $\sigma$ contains a
subsequence isomorphic to $\tau$. Otherwise, we say that $\sigma$
\emph{avoids} $\tau$. The \emph{reversal} of $\tau$, denoted by $r(\tau$), is
the pattern $\tau$ read from right to left, and the
\emph{complement} of $\tau$, denoted by $c(\tau)$, is the pattern
obtained by replacing $\tau_i$ by $k+1-\tau_i$. The set
$\{\tau,r(\tau),c(\tau),c(r(\tau))\}$ is called the
\emph{symmetry} class of $\tau$. It is easy to see that patterns
from the same symmetry class occur an equal number of times in all
the words of length $m$.

The connection between compositions and words is as follows: If $C_\tau^A(x,y,z)$ is the generating function for the
number of compositions of $n$ with $m$ parts in the set $A$ and
$r$ occurrences of the pattern or statistic $\tau$, then
$C_\tau^A(1,y,z)$ is the generating function for the number of
words of length $m$ on the alphabet $A$ with $r$ occurrences of
the subword or statistic $\tau$.  We are now ready to apply our results from the previous sections to words on $k$ letters. 

Theorem~\ref{th111} gives
    $$C_{111}^{[k]}(1,y,z)=\frac{1}{1-\frac{kz(1+(1-y)z)}{1+z(1+z)(1-y)}}=
    \frac{1+z(1+z)(1-y)}{1-(k-1+y)z-(k-1)(1-y)z^2},$$
i.e., we obtain the results of Example 2.2 \cite{BurMan2003a} and Theorem 3.1 \cite{BurMan2003}.

Theorem~\ref{th112} gives (after simplification)
$$C_{112}^{[k]}(1,y,z)=C_{221}^{[k]}(1,y,z)=
\frac{(1-y)z}{(1-y)z-1+(1-(1-y)z^2)^k},$$ i.e., we obtain Theorem
3.2 \cite{BurMan2003}. In addition, for $y=0$ we obtain the generating
function for the words of length $m$ that avoid the pattern $112$,
given in Theorem 3.10 \cite{BurMan2003a}:
$$C_{112}^{[k]}(1,0,z)=C_{221}^{[k]}(1,0,z)=\frac{z}{z-1+(1-z^2)^k}.$$

Note that necessarily $C_{221}^{[k]}(1,y,z)= C_{112}^{[k]}(1,y,z)$, 
as 112 and $221=c(112)$ are in the same symmetry class, which explains why the respective generating functions for compositions are similar in structure.

Theorem~\ref{th123} for $x=1$ yields $t_k^p([k])=\binom{k}{p}z^p$
and thus,
$$C_{123}^{[k]}(1,y,z)=\frac{1}{1-k\,z-\sum_{p=3}^k\sum_{j=0}^{p-3}\binom{p-3}{j}\binom{k}{p+j}z^{p+j}(y-1)^{p-2}}.$$
This generating function for the number of words of length $m$
that contain the pattern $123$ exactly $r$ times was given in a
different form in Theorem 3.3 of~\cite{BurMan2003}:

{\bf Theorem 3.3}~\cite{BurMan2003} {\emph For $k \ge 2$,
$$F_{123}(z,y;k)= \frac{1}{1-k\,z-\sum_{j=3}^k(-z)^j\binom{k}{j}(1-y)^{\lfloor j/2 \rfloor}U_{j-3}(y)},$$
where $U_0(y)=U_1(y)=1,U_{2n}(y)=(1-y)U_{2n-1}(y)-U_{2n-2}(y)$,
and $U_{2n+1}(y)=U_{2n}(y)-U_{2n-1}(y)$. Furthermore, the generating function
for $U_n(y)$ is given by
$$\sum_{n \ge 0}U_n(y)z^n=\frac{1+z+z^2}{1+(1+y)z^2+z^4}.$$}

We can prove the equivalence of the two generating functions by
substituting $x=1$ into equations (\ref{CD123}) and (\ref{DA1})
which yields the expressions given in \cite{BurMan2003} for
$F_{\tau}(z,y;k)$ and $D_{\tau}(z,y;k)$. Comparison of the initial
conditions then shows that $C_{123}^{[k]}(1,y,z)=F_{123}(z,y;k)$.

Substituting $y=0$ and $x=1$ in Theorem~\ref{th123} we get the
generating function for the number of words of length $m$ avoiding
the subword $123$:
$$C_{123}^{[k]}(1,0,z)=\frac{1}{1-k\,z-(-1)^p\sum_{p=3}^k\sum_{j=0}^{p-3}\binom{p-3}{j}\binom{k}{p+j}z^{p+j}}.$$

This can be shown to give the result of Theorem 3.13~\cite{BurMan2003a},
namely
\begin{equation}\label{F123}
 F_{123}(z,0;k)= \frac{1}{\sum_{j=0}^k a_j \binom{k}{j} z^j},
 \end{equation}
where $a_{3\ell} = 1$, $a_{3\ell+1} = -1$ and $a_{3\ell+2}=0$. We use the form of Theorem 3.3~\cite{BurMan2003} for $y=0$, which gives
$$F_{123}(z,0;k)= \frac{1}{1-k\,z+\sum_{j=3}^k(-1)^{j+1}z^j\binom{k}{j}U_{j-3}(0)}.$$
Note that for $y=0$ and $m \ge 0$, $U_m(0)=U_{m-1}(0)-U_{m-2}(0)$, and thus, $U_m(0)=1$ for $m \equiv 0$ or $1$(mod 6), $U_m(0)=-1$ for $m \equiv 3$ or $4$(mod 6), and $U_m(0)=0$ otherwise.
For $j\ge 3$ and $j=3m$, $U_{j-3}(0)=1$ and $(-1)^{j+1}=1$. If $j \ge 3$ and $j=3m+1$, then $U_{j-3}(0)$ and $(-1)^{j+1}$ will have opposite signs, resulting in a coefficient of $-1$. Finally, if $j \ge 3$ and $j=3m+2$, then $U_{j-3}(0)=0$. Note that formula~(\ref{F123}) also holds for $j=0,1,2$, and therefore we have shown the equivalence of the two results.

Now we apply our results to the statistics $peak$ and
$valley$, which will give new results. Using Theorem~\ref{peak}
for $A=[k]$ and $x=1$, we need to determine the products
$M^s([k])=\sum_{(i_1,i_2,\ldots,i_s) \in P^s ([k])} \prod_{j=1}^s
b_{i_j}=z^s|P^s([k])|$. To determine $|P^s([k])|$ for $s=2\ell+1$ and $s = 2\ell$,
we rewrite the sequence of alternating $\le$ and $<$ signs in the
definition of $P^s([k])$ as strict inequalities and obtain (with $j_{n}=i_{n}+\lfloor{(n-1)/2}\rfloor$)
$$|P^{2\ell+1}([k])|=
|\{(j_1,j_2,\ldots,j_{2\ell+1})| 1\leq
j_1<j_2<\cdots<j_{2\ell+1}\leq k+\ell\}|=\binom{k+\ell}{2\ell+1}.$$
 Using the same argument, we obtain
$|P^{2\ell}([k])|=
\binom{k-1+\ell}{2\ell}$.
Substituting $M^{2\ell}([k])=z^{2\ell}\binom{k-1+\ell}{2\ell}$ and
$M^{2\ell+1}([k])=z^{2\ell+1}\binom{k+\ell}{2\ell+1}$ into Theorem~\ref{peak} gives
\begin{align*}
C_{peak}^{[k]}(1,y,z)
&=\frac{\sum_{j\geq0}z^{2j}(1-y)^j\binom{k-1+j}{2j}}
{\sum_{j\geq0}z^{2j}(1-y)^j\binom{k-1+j}{2j}-\sum_{j\geq0}z^{2j+1}(1-y)^j\binom{k+j}{2j+1}}.
\end{align*}
Since $c(peak)=valley$, we have $C_{valley}^{[k]}(1,y,z)=C_{peak}^{[k]}(1,y,z)$, and, setting $y=0$, the generating function for the number of words of length
$m$ on the alphabet $[k]$ without peaks (valleys) is given by
$$C_{peak}^{[k]}(1,0,z)=C_{valley}^{[k]}(1,0,z)=\frac{\sum_{j\geq0}z^{2j}\binom{k-1+j}{2j}}
{\sum_{j\geq0}z^{2j}\binom{k-1+j}{2j}-\sum_{j\geq0}z^{2j+1}\binom{k+j}{2j+1}}.$$

Note that once again the symmetry structure of words explains the fact that the generating functions for the number of compositions that avoid $valleys$ and $peaks$, respectively, have similar structure.

\vskip 30pt

\section*{Acknowledgments}
 The authors would like to thank the anonymous referee for his/her careful reading of the manuscript and for constructive suggestions that have improved this paper, including references on Carlitz compositions using a probabilistic approach~\cite{GohHit2002, HitLou2001, LouPro2002}. 
\vskip 30pt

\end{document}